\documentclass[letterpaper, 10 pt, conference]{ieeeconf}  

\IEEEoverridecommandlockouts                              

\usepackage{graphics} 
\usepackage{epsfig} 
\usepackage{mathptmx} 
\usepackage{times} 
\usepackage{amsmath} 
\usepackage{amssymb}  
\usepackage{multirow}
\usepackage{url}

\usepackage{cite}

\title{\LARGE \bf
Modeling, scientific computing and optimal control \\ for renewable energy systems with storage
}

\author{Nicola Cantisani, Tobias K. S. Ritschel, Christian A. Thilker,
Henrik Madsen and John Bagterp J{\o}rgensen
\thanks{N. Cantisani, T.K.S. Ritschel, C.A. Thilker, H. Madsen, J.B. Jørgensen are with Department of Applied Mathematics and Computer Science, Technical University of Denmark, DK-2800 Kgs. Lyngby, Denmark}%
}

\begin{document}

\maketitle
\thispagestyle{empty}
\pagestyle{empty}

\begin{abstract}
This paper presents models for renewable energy systems with storage, and considers its optimal operation. We model and simulate wind and solar power production using stochastic differential equations as well as storage of the produced power using batteries, thermal storage, and water electrolysis. We formulate an economic optimal control problem, with the scope of controlling the system in the most efficient way, while satisfying the power demand from the electric grid. Deploying multiple storage systems allows flexibility and higher reliability of the renewable energy system.
\end{abstract}

\section{Introduction}
\label{sec:Introduction}
\begin{figure}[hb]
    \centering
    \includegraphics[width=0.5\textwidth]{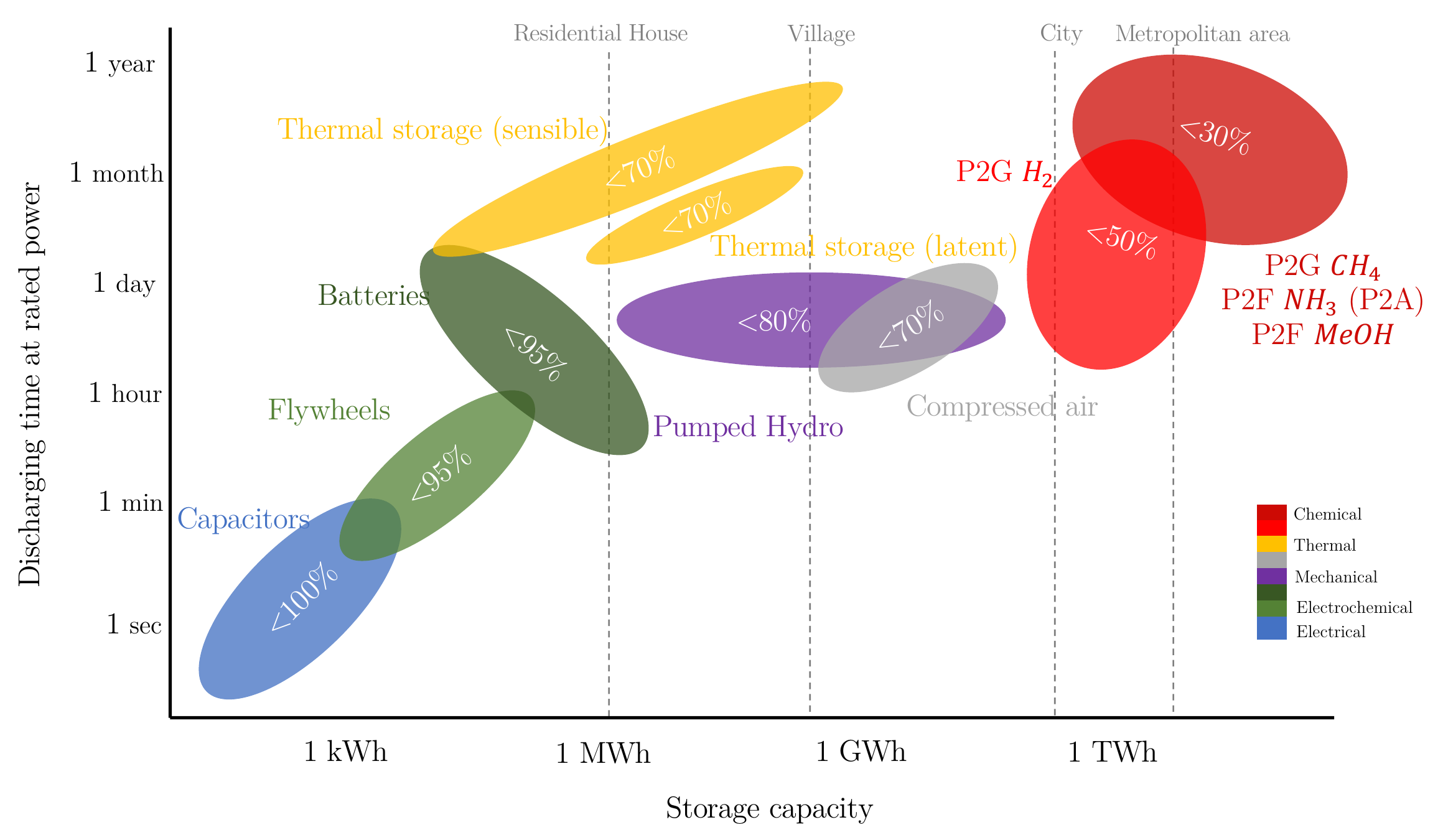}
    \caption{Characterization of the main energy storage systems (ESS) by storage capacity and discharging time at rated power (i.e. maximum output power). The percentages refer to energy efficiency.}
    \label{fig:energy_storage_systems}
\end{figure}

As society has progressively become industrialized and the use of energy has been increasing over the last century, greenhouse gas emissions have risen at an alarming pace, leading to global warming. The biggest contribution to greenhouse emissions is carbon dioxide, which is mainly released when burning fossil fuels. Consequently, we need to push development towards the use of renewable energy sources (RESs), such as wind and solar. However, solar and wind cannot be controlled and are highly volatile/intermittent. Energy storage systems (ESSs) can be a solution to this problem, by serving as a buffer when the power production is low or absent, and making the energy more dispatchable \cite{Braga:etal:2021,Krishan:etal:2019}. Fig. \ref{fig:energy_storage_systems} presents the main ESSs, providing information about the capacity, the discharging time at rated power, and the energy conversion efficiency. Even though ESSs with bigger capacities generally have lower energy efficiencies, they can support energy systems for longer periods of time. 

In this paper, we demonstrate mathematical modeling, simulation, and control of such systems. Based on the scientific literature, we model the uncertainty coming from the wind and the sun using stochastic differential equations. We express the power extraction from wind turbines as an optimization problem and present models for energy storage. We develop and propose simple models for battery and thermal storage, based on energy balance equations. Electrolysis is modeled using an electrochemical model, as in \cite{SAKAS20224328,ULLEBERG200321}. We will consider a hybrid energy system consisting of stochastic power production from wind and solar and the combination of battery, thermal storage, and hydrogen for energy storage.

\begin{figure}[b]
    \centering
    \includegraphics[width=0.49\textwidth]{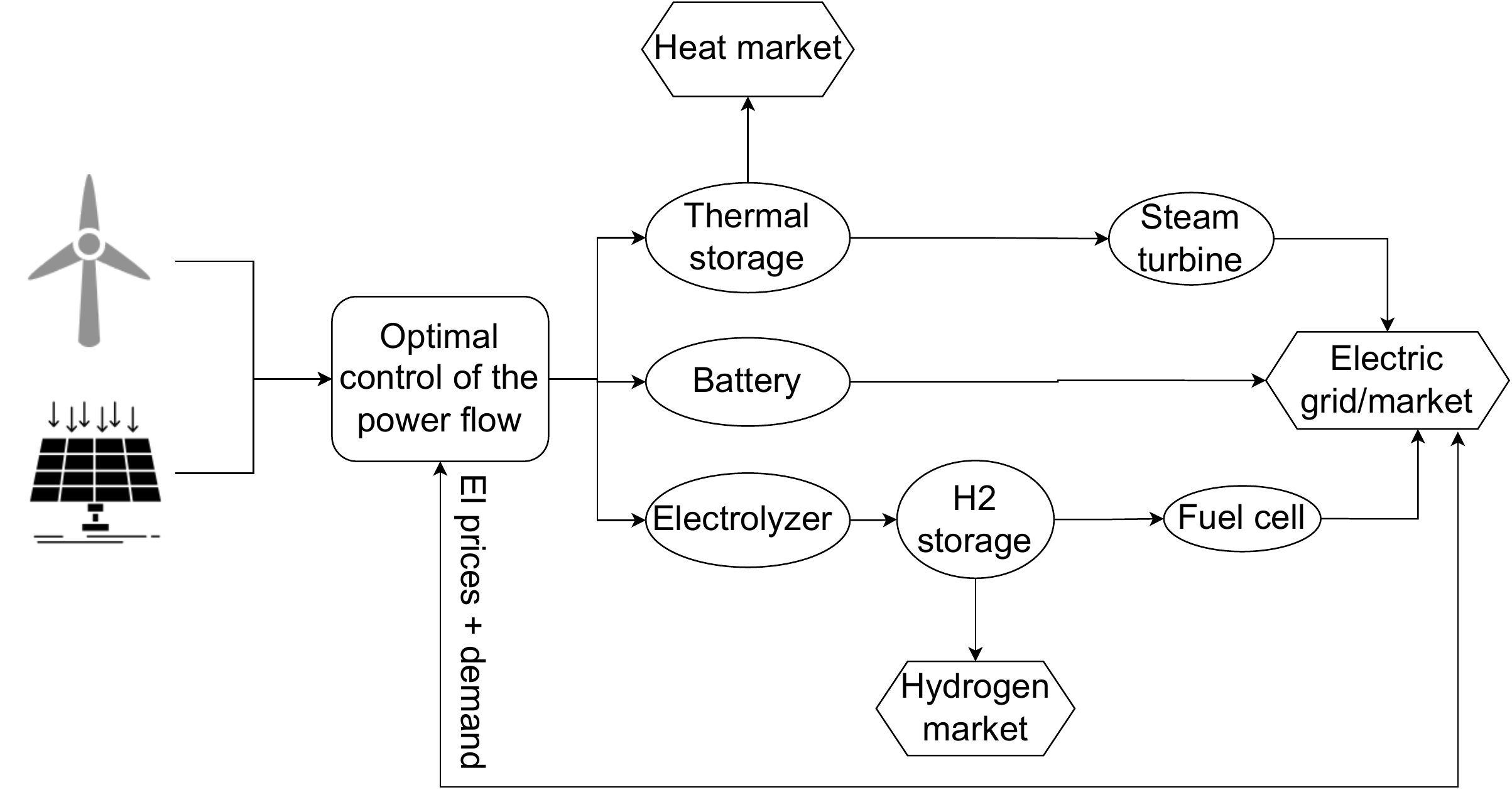}
    \caption{Overview of the controllable renewable energy system.}
    \label{fig:diag}
\end{figure}

Ultimately, we want to achieve optimal operation of such a system. Therefore, an optimal control problem is formulated. The goal is to maximize the profit generated from selling electricity, heat, and hydrogen, while assuring that the grid's power demand is always satisfied. Fig. \ref{fig:diag} provides an overview of the entire energy system and the controllable power flows. 

There is growing research on how to best control a renewable energy system with ESSs, where model predictive control (MPC) is often the preferred choice. \cite{Belloni:etal:2016,Khatamianfar:etal:2012,Malkowski:etal:2018,Teleke:etal:2009} show how to control a system with one stochastic power source (wind or solar) and a single storage device, often a battery. Integrating and controlling different storage solutions together is also attempted in the literature \cite{Trifkovic:etal:2013,Aldaouab:etal:2018}. The control strategy proposed by the mentioned resources is based on maximizing the utilization of renewable resources. In this work, we explore the possibility of integrating and controlling three different ESSs, maximizing the economic revenue when such a system is interfaced with the market.

The paper is structured as follows. Sections \ref{sec:WindEnergy} and \ref{sec:SolarEnergy} present models for power production from wind turbines and solar photovoltaic (PV) panels, respectively. In Sections \ref{sec:Battery}, \ref{sec:ThermalStorage} and \ref{sec:Electrolysis} we model energy storage via battery, thermal storage and water electrolysis, respectively. Section \ref{sec:optimal_control} introduces the optimal control problem. Finally, Section \ref{sec:simulation} presents a simulation of the system.

\section{Wind energy} \label{sec:WindEnergy}
In this section, we present control relevant simulation models for horizontal axis wind turbines (HAWT). The model consists of a model for the wind speed and one for the wind turbine.

\subsection{Wind speed model}
The wind speed is modeled as the sum of two components \cite{Soltani:etal:2013}:
\begin{equation}
	v(t) = v_m(t) + v_t(t).
\end{equation}
$v_m(t)$ is the mean wind speed and ${v}_t(t)$ is the turbulent wind speed. Each part is modeled by a stochastic differential equation (SDE):
\begin{subequations}\label{wind_speed_model}
\begin{align}
	d v_t(t) & = - \frac{\pi \, v_m(t)}{2L} v_t(t) \, dt +  \sqrt{ \frac{\pi \, v_m(t)^3 \, t_i^2}{L}} \,  d w_1(t),\\
    d v_m (t) & = \sigma_2 \, dw_2(t).
\end{align}
\end{subequations}
$L$ is the turbulence length and $t_i$ is the turbulence intensity. The parameter values for the wind speed model are reported in Table \ref{tab:wind_speed}. Notice that the mean wind speed can also be set and fixed for simulation purposes. Observe that we do not consider the wind direction here, but it is relevant for a wind farm. 
\begin{table}[tb]
\centering
\caption{Parameters of the wind speed model \cite{Soltani:etal:2013}.}
    \begin{tabular}{|c|l l|}
        \hline
        \textbf{Parameter} & \textbf{Value} & \textbf{Unit} \\
        \hline
	$L$ & 170.1 & m \\
	 $t_i$ & 0.2 & - \\
		 $\sigma_2$ & $\sqrt{4/600}$ & m $\cdot$ s$^{-3/2}$ \\
		\hline
	\end{tabular}
	\label{tab:wind_speed}
\end{table}

\subsection{Power extraction}
The available power in the wind is
\begin{equation}
	P_w(t) = \frac{1}{2} \rho \pi R^2 v(t)^3.
\end{equation}
$\rho$ is the density of the air and $R$ is radius of the blade. The tip speed ratio is the ratio between the tangential speed at the tip of a blade of the turbine and the wind speed, i.e.
\begin{equation}
\lambda(t) = \frac{R \Omega_r(t)}{v(t)}.
\end{equation}
$\Omega_r(t)$ is the rotational speed of the rotor.
The blade pitch angle, $\theta$, is defined as the angle between the chord of the airfoil and the rotational plane.
The extracted power by the wind turbine rotor is only a part of the full available power, $P_w(t)$. The fraction of the extracted power depends on the turbine's power coefficient, $C_P(\lambda(t),\theta(t))$, i.e.
\begin{equation}
	P_r(t) = P_w(t) C_P(\lambda(t),\theta(t)).
\end{equation}
The power coefficient values are table-based.
The torque is 
\begin{equation}
	Q_r(t) = \frac{P_r(t)}{\Omega_r(t)}.
\end{equation}

\subsection{Stationary optimal operation}
\begin{table}[tb]
	\centering
	\caption{Relevant NREL 5MW parameters and constraints \cite{nrel}.}
    \begin{tabular}{|c|l l|l|} \hline
    \textbf{Parameter} & \textbf{Value} & \textbf{Unit} & \textbf{Description} \\
    \hline
	$\rho$ & 1.225 & kg/m$^3$ & Density of air \\
	$R$ & 62.94 & m & Length of rotor blade	\\ 
	$\eta_g$ & 94.4\% & - & Power generator efficiency \\ 
	$v_{cut-in}$ & 3 & m/s & Cut in velocity \\
	$v_{rated}$ & 11.4 & m/s & Rated velocity \\
	$v_{cut-out}$ & 25 & m/s & Cut out velocity \\
	$\Omega_{r,min}$ & 6.9 & rpm & Minimum rotor speed \\
	$\Omega_{r,max}$ & 12.1 & rpm & Maximum rotor speed  \\
	$\theta_{min}$ & -5 & deg & Minimum pitch angle  \\
	$\theta_{max}$ & 25 & deg & Maximum pitch angle \\ \hline
    \end{tabular}	
\label{tab:5mwpar}
\end{table}

The maximum power extraction of the turbine is achieved as the solution of the following optimization problem, assuming steady state conditions (i.e. constant wind speed).
\begin{subequations}
\begin{alignat}{3}
	& \max_{\theta,\Omega_r} \quad && \frac{1}{2} \rho R^2 v^3 C_P(\lambda(\Omega_r,v),\theta) \\
	& \mathrm{s.t.} && \theta_{\min} \leq \theta \leq \theta_{\max}, \label{eq:constraint_pitch}\\
	&       && \Omega_{r,\min} \leq \Omega_r \leq \Omega_{r,\max}. \label{eq:constraint_omega}
\end{alignat}
\end{subequations}
Notice that the problem depends on the wind speed, $v$. Every value of the wind speed will lead to a different solution. The solution for every wind speed between the cut-in and cut-out velocity produces the power curve of the wind turbine.
The actual generated power by the generator, $P_g$, in the turbine is obtained by multiplying the power produced by the rotor, $P_r$, by the efficiency of the electrical generator, $\eta_g$.
The turbine is, in general, constrained between a range of pitch angles and rotational speed, according to the manufacturer's regulations, hence constraints \eqref{eq:constraint_pitch} and \eqref{eq:constraint_omega}.\\
We can reformulate the optimization problem, emphasizing the constraint using $\lambda$, and discarding the constant part of the objective function.
\begin{subequations}\label{optimization_problem}
	\begin{alignat}{3}
	& \max_{\theta,\lambda} \quad && C_P(\lambda,\theta) \\
	& s.t. && \theta_{\min} \leq \theta \leq \theta_{\max}, \\
	&       && \lambda_{\min}(v) \leq \lambda \leq \lambda_{\max}(v),
	\end{alignat}
\end{subequations}
where
\begin{equation}
    \lambda_{\min}(v) = \frac{R \Omega_{r,\min}}{v}, \quad	\lambda_{\max}(v) = \frac{R \Omega_{r,\max}}{v}.
\end{equation}

The solution, for each wind speed $v$, consists of $\theta^* = \theta^*(v)$, $\lambda^* = \lambda^*(v)$, and the derived quantities $C_P^*(v), \, P_w^*(v), \, P_g^*(v)$ and $\Omega_r^*(v).$

\begin{figure}[tb]
    \centering
    \includegraphics[width=0.5\textwidth]{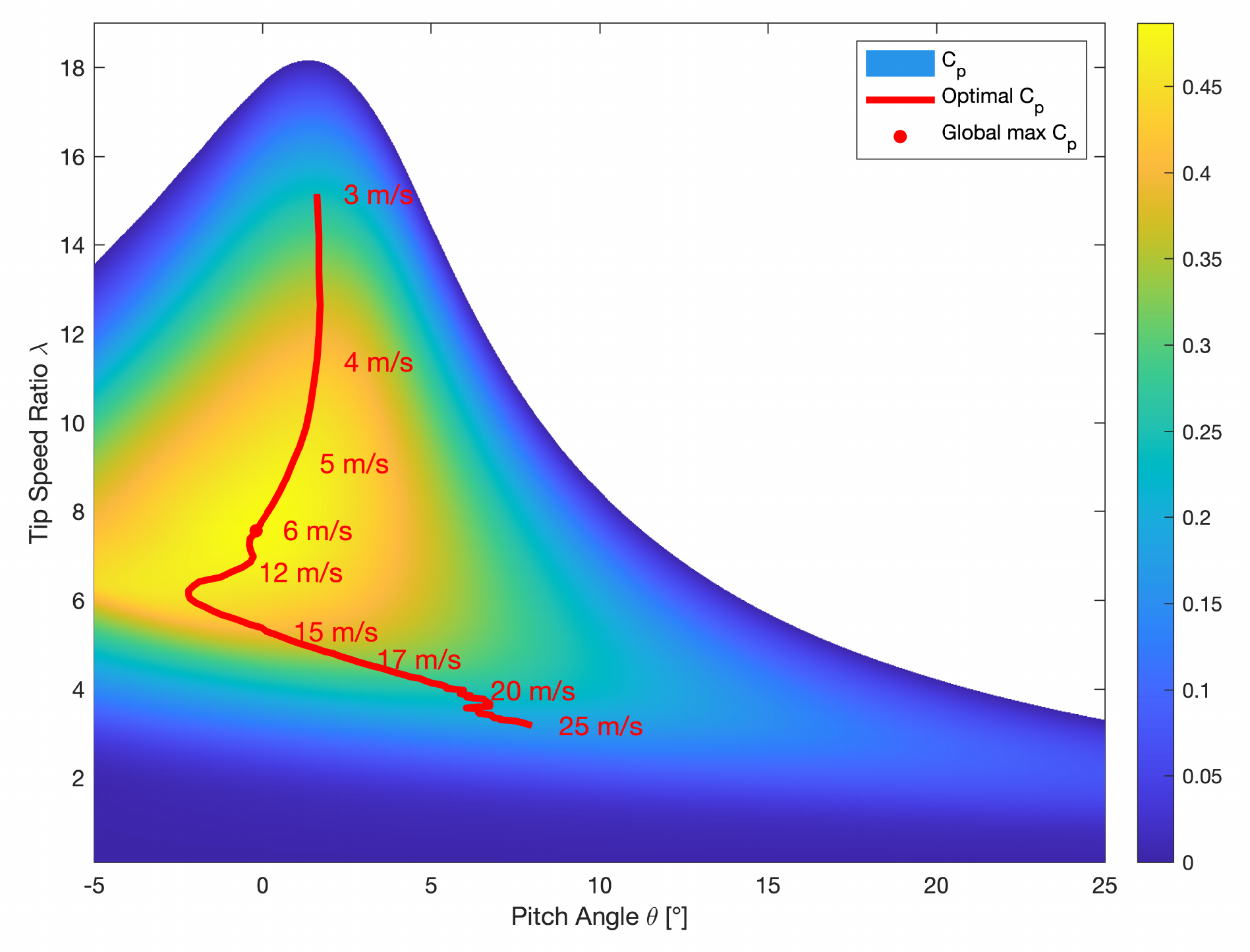}
    \caption{$C_p(\lambda,\theta)$ values of a NREL 5MW wind turbine and the solution of the optimal stationary problem \ref{optimization_problem}, for $v \in [v_{cut-in},v_{cut-out}]$.}
    \label{fig:cp}
\end{figure}
\begin{figure}[tb]
    \centering
    \includegraphics[width=0.5\textwidth]{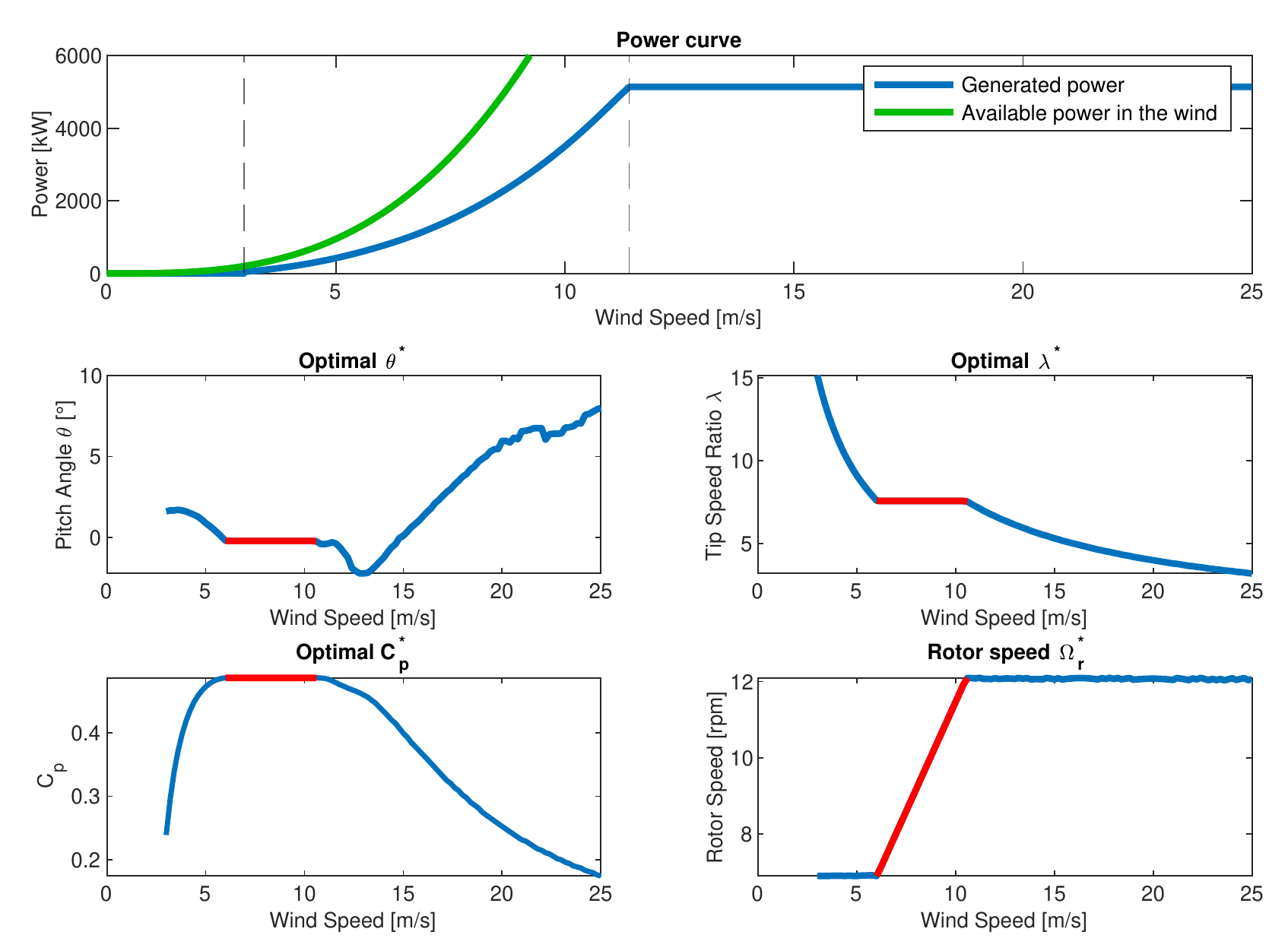}
    \caption{The optimal pitch angle $\theta^*$, tip speed ratio $\lambda^*$, power coefficient $C_p^*$, rotor speed $\Omega_r^*$ and generated power $P_g^*$ as a function of the wind speed.}
    \label{fig:sol}
\end{figure}

We consider an NREL 5MW wind turbine \cite{nrel} and solve the problem, using the parameters in Table \ref{tab:5mwpar}. The $C_p(\lambda,\theta)$ values and the solution of the optimal stationary problem \eqref{optimization_problem} for every $v$ is shown in Fig. \ref{fig:cp}. Fig. \ref{fig:sol} reports the entire solution as a function of the wind speed.

\section{Solar Energy} \label{sec:SolarEnergy}

The global radiation, $I_G$, is the total amount of short-wave radiation received from the sun by a surface horizontal near the ground. This value is of particular interest to photovoltaic installations and depends on both direct radiation, $I_N$, and diffuse radiation, $I_D$. 

In the literature, the solar radiation is commonly calculated as the intensity that would be obtained in clear sky conditions, multiplied by a factor that depends on the cloud cover \cite{paltridge1976a,madsen1985a,Thilker:etal:2021,taesler1984a,bacher2009a}. Consequently, simulation models for direct, diffuse, and global solar radiation have been suggested in the literature where stochastic models for cloud cover have been used as the driving entity \cite{madsen_phd,madsen1985a,Thilker:etal:2021}. Such methods have been used for simulation and prediction of weather observations in optimal control of heat supply to buildings \cite{madsen87a}.

In this paper, the cloud cover $\kappa(t)$ is modeled using the mean-reversion SDE model from \cite{Thilker:etal:2021}
\begin{subequations}
\begin{align}
\begin{split}
    d\kappa(t) &= \theta(\kappa(t)) [\mu(\kappa(t))-\kappa(t) ] \, dt \\ & \qquad + \sigma \kappa(t)[1-\kappa(t)] \, d\omega_{\kappa}(t),
\end{split} \\
    \mu(\kappa(t)) &=\frac{\exp[P_7(\kappa(t))]}{1+\exp[P_7(\kappa(t))]}, \\
    \theta(\kappa(t)) &= \tilde{\theta} \sqrt{\kappa(t) (1-\kappa(t))}.
\end{align}
\end{subequations}
The term $P_7(.) = \sum_{k=1}^7 p_k L_k(.)$ is the linear combination of the first seven Legendre polynomials. The equation provides a cloud cover such that $\kappa(t) \in [0,1]$. A transformation of the data is needed to map the data to the Okta scale ($0$ = absence of cloud cover, $8$ = total coverage of the sky). For simulation purposes, the mean value, $\mu$, can also be fixed.

We now present the dynamical models that describe the solar irradiance \cite{madsen_phd}.
The model for the direct radiation, $I_N$, is 
\begin{subequations}
\begin{align}
    \frac{I_N(\kappa(t),h(t))}{I_{0,N}(h(t))} & = f_N(\kappa(t)) + \epsilon_N(t), \\
    I_{0,N}(h(t)) & = a_N [1-\exp(-b_N h(t))], \\
    f_N(\kappa(t))& = \frac{a'+b' \cos(\frac{2 \pi d_t}{365} + c')}{1+\exp(\kappa(t)-\alpha_N)},\\
    \epsilon_N(t) & \sim \mathcal{N} \left( 0, r_N^2\frac{\sigma_{D,\kappa(t)}^2}{\sin^2h(t) I_{0,N}(h(t))^2} dt \right).
\end{align}
\label{eq:direct_radiation}
\end{subequations}
$I_{0,N}$ is the direct radiation in clear skies (non-stochastic), $d_t$ is the day of the year and $h(t)$ is the sun elevation angle at the location of interest. The cloud cover is here provided in Okta scale ($\kappa(t) \in [0,8]$). More variability is expected when the sun elevation, $h(t)$, is low. Madsen and Thyregod \cite{madsen1988a} suggest a stochastic model for direct radiation in clear skies which also takes into account the auto-correlation caused by the persistence in variations of the transparency of the atmosphere due to ozone, dust, water, etc.

The model for the diffuse radiation, $I_D$, is
\begin{subequations}
\begin{align}
    \frac{I_D(\kappa(t),h(t))}{I_{0,D}(h(t))} & = f_D(\kappa(t)) + \epsilon_D(t), \\
    I_{0,D}(h(t)) & = c_D + a_D [1-\exp(-b_D h(t))], \\
    \begin{split}
    f_D(\kappa(t))& = a_0(t)+ a_1(t)(1-\kappa(t)/8) \\ & \qquad +
    a_2(\kappa(t)/8)^{\alpha_D}(1-\kappa(t)/8),
    \end{split} \\
    \epsilon_D(t) & \sim \mathcal{N} (0, \sigma_{N,\kappa(t)}^2 r_D^2 \, dt).
\end{align}
\end{subequations}
$I_{0,D}$ is the diffuse radiation in clear skies and
\begin{subequations}
\begin{align}
    a_0 = r_1 + r_2  \cos \left( \frac{2 \pi d_t}{365} + r_3\right),\\
    a_1 = k_1 + k_2  \cos \left( \frac{2 \pi d_t}{365} + k_3 \right).
\end{align}
\end{subequations}
The variances, $\sigma_{N,\kappa(t)}$ and $\sigma_{D,\kappa(t)}$, depend on the value of the cloud cover, $\kappa(t)$ (data is available in \cite{madsen_phd}).

The global radiation is
\begin{equation}
    I_G(t) = I_N(\kappa(t),h(t)) \sin h(t) + I_D(\kappa(t),h(t)).
\end{equation}
Notice that, due to the dependency of the solar elevation in the variance structure of \eqref{eq:direct_radiation}, the variance of the global radiation becomes independent of the solar elevation. 

The power generated by solar radiation is
\begin{equation}
    P_s(t) = \eta_s A_s I_G(t),
\end{equation}
where $A_s$ is the area of PVs and $\eta_s$ the PVs efficiency. For real-time solar power forecasting tools, like SOLARFOR\footnote{\url{https://enfor.dk/services/solarfor/}}, the efficiency depends, e.g., on the dirtiness of the PV panels. 
Table \ref{tab:solar_energy} reports all the parameters of the models.
\begin{table}[tb]
\centering
\caption{Parameters of the solar power model.}
    \begin{tabular}{|llllll|}
        \hline
        \multicolumn{6}{|c|}{\textbf{Cloud cover model} \cite{Thilker:etal:2021}}\\
        \hline
        Parameter    & $\tilde{\theta}$  & $\sigma$ & $p_1$ & $p_2$ &$p_3$ \\
        Value        & 0.187 & 0.835 & -53.1 & 14.6 & -42.3  \\
 		\hline
        Parameter            & $p_4$ &$p_5$ & $p_6$ & $p_7$ &  \\
        Value            &  8.8  & -58.1 & -30.3 & -45.7 & \\
 		\hline
 		 \multicolumn{6}{|c|}{\textbf{Direct radiation model} \cite{madsen_phd}}\\
         \hline
        Parameter    & $a_N$ & $b_N$ & $r_N$ & $\alpha_N$ & $a'$ \\       
        Value        &  842.3  & 0.0614  & 1.0430 & 4.6368 & 1.1354 \\
        \hline
        Parameter             & $b'$ & $c'$ & & & \\
        Value             & 0.1965  & -0.2571 & & & \\
        \hline
        \multicolumn{6}{|c|}{\textbf{Diffuse radiation model} \cite{madsen_phd}}\\
         \hline
        Parameter    & $a_D$ & $b_D$ &$c_D$ & $r_1$ & $r_2$ \\
        Value        & 161.1 & 0.0333 & 3.68 & 0.7067 & -0.2456  \\
        \hline       
        Parameter             & $r_3$ & $k_1$ & $k_2$  &$k_3$ & $a_2$ \\
        Value             & 0.5625 & 0.1946 & 0.1549 & 0.6034 & 6.7033 \\
        \hline
        Parameter             & $\alpha_D$ & $r_D$ & & & \\
        Value             &  2.2993 & 1.0170 & & & \\
        \hline
	\end{tabular}	
	\label{tab:solar_energy}
	
\end{table}

  
\section{Battery} \label{sec:Battery}


We model the accumulated energy in a battery considering its energy balance. The change in the energy in the battery $E_b(t)$ over time is expressed by the following differential equation
\begin{equation} \label{eq:battery}
\dot E_{b}(t) = -\alpha_b E_b(t) + \eta_{b,in} P_{b,in}(t) - \frac{1}{\eta_{b,out}} P_{b,out}(t),\\
\end{equation}
where $P_{b,in}(t)$ is the incoming power and $P_{b,out}(t)$ is the outgoing power. The factors $\eta_{b,in}$ and $\eta_{b,out}$ are the charging and discharging efficiencies. The first term of the equation accounts for the loss of charge over time. Notice that the charge is constrained between a minimum and a maximum value,
\begin{equation}
    E_{b,min} \leq E_b(t) \leq E_{b,max}.
\end{equation}
\section{Thermal storage} \label{sec:ThermalStorage}
Thermal energy storage (TES) is a relatively inexpensive and efficient way to store big amounts of energy for a long period of time. Even though the TES systems can differ quite significantly, they all use the ability of some material to accumulate energy through heat exchange. When energy is needed, the material releases heat (cools down), which can be used in several ways. 

TES can be either latent or sensible, according to whether a phase transition of the medium is allowed or not. We consider the latter case. In this case, a heat exchange with the material always implies a change in temperature. We express the heat exchange using the first law of thermodynamics (in absence of work):
\begin{equation}
    \Delta U = Q = m \, c_P \, \Delta T.
\end{equation}
$m$ is the mass and $c_P$ is the specific heat capacity. Common materials used for TES are molten salts. \cite{Caraballo} reports density, fusion temperature, $T_{\min}$, decomposition temperature, $T_{\max}$, and specific heat capacity, $c_P(T)$, for the most commonly used molten salts. We can compute how much energy it is possible to store in a unit of mass of the material, considering the specific storage capacity (J/kg), as
\begin{equation}
    \Delta u = \int_{T_{\min}}^{T_{\max}} c_p(T) dT.
\end{equation}
As for the battery, the thermal storage model is based on energy balance. The factors contributing to an internal energy change are incoming power, $P_{t,in}$, heat flow (loss) to the ambient, $Q_{ta}$, and outgoing heat flow, $Q_{t,out}$.
\begin{equation}
	\dot U(t) = \eta_{t,in} P_{t,in}(t) - Q_{ta}(t)  - Q_{t,out}(t).
\end{equation}
Expressing all the terms, we get the full model
\begin{equation} \label{eq:thermal_storage}
\begin{split}
& m \, c_P(T(t)) \dot T(t) = \\ & \quad \eta_{t,in} P_{t,in}(t) - (UA)_t (T(t)-T_a(t)) - Q_{t,out}(t), 
\end{split}
\end{equation}
subject to 
\begin{equation}
    T_{\min} \leq T(t) \leq T_{\max}.
\end{equation}
$(UA)_t$ is the overall conductivity of the storage tank and depends on the material and the exchange area. $T_a(t)$ is the ambient temperature.

\section{Water electrolysis}\label{sec:Electrolysis}
\subsection{Thermodynamics}
The general reaction of electrolysis of water is
\begin{equation}
    \mathrm{2H_2O (l) + \text{Electrical power} \longrightarrow 2H_2 (g) +O_2 (g)}.
\end{equation}
The thermodynamics of the reaction needs to be considered. The energy that is required for the reaction to happen is the enthalpy of the process
\begin{equation}
    \Delta H = \Delta G + T \Delta S.
\end{equation}
The reaction is endothermic because $\Delta H >0$ and it is equal to $\Delta H^0 = 285.8$ kJ/mol at standard conditions. The change in Gibbs free energy is equal to $\Delta G^0 = 237.2 $ kJ/mol at standard conditions. A positive value in the change of Gibbs free energy indicates that the reaction does not happen spontaneously, but it is spontaneous in the opposite direction. To let the reaction happen in this direction, electric work must be supplied, and it is equal to the change in Gibbs free energy. The remaining energy requirement is brought as heat, because of the irreversibilities of the process and it is equal to $Q = T \Delta S = 48$ kJ/mol. We can re-write the total reaction, including the energy demands
\begin{equation}
\begin{split}
    \mathrm{ H_2 O (l) + 237.2\, \text{kJ}_{\text{electricity}}/\text{mol} + 48.6\,\text{kJ}_{\text{heat}} /\text{mol}} \\\mathrm{  \longrightarrow H_2 (g) + \frac{1}{2} O_2 (g)}.
\end{split}
\end{equation}

The lowest voltage that is needed to initiate the process is called reversible voltage. It is given by
\begin{equation}
    U_{rev} = \frac{\Delta G}{z F}.
\end{equation}
$z=2$ is the number of electrons transferred and $F = 96485$ C/mol is the Faraday constant. At standard conditions, $U_{rev}^0 = 1.229$ V.

In order for the reaction to happen, heat must also be provided. In real electrolyzers, additional electric work is provided, instead of heat. In this case, the minimum voltage requirement is higher, and it is called the thermoneutral voltage
\begin{equation} \label{eq:tn}
    U_{tn} = \frac{\Delta H}{z F}.
\end{equation}
At standard conditions, $U_{tn}^0 = 1.481$ V. Both the reversible and thermoneutral voltages are state functions.

\subsection{Electrochemical model}
We consider the electrochemical model from \cite{ULLEBERG200321}. This has been used extensively in the scientific literature, because of its simplicity and its empiric basis. We consider only alkaline electrolyzers.
We model the kinetics at the electrodes, i.e. we investigate the relationship between current intensity and voltage in the electrolyzer ($I-U$).
In general, we can express the \textit{real} total voltage in one electrolytic cell as the sum of several terms
\begin{equation} \label{eq:el1}
    U_{cell} = U_{rev} +  U_{ohm} + U_{act} + U_{con}.
\end{equation}
$U_{act}$ is the activation overvoltage, which is caused by the anode and cathode half-reactions kinetics \cite{Koponen}. $U_{ohm}$ is the ohmic overvoltage. This is due to ohmic losses in the cell, generated by the resistance to the electron flux in the components of the electrolyzer. $U_{con}$ is concentration overvoltage, which is negligible for alkaline electrolyzers.
\cite{ULLEBERG200321} proposed the following electrochemical model for ohmic and activation overvoltages
\begin{subequations} \label{eq:el2}
    \begin{align}
    U_{ohm} & = (r_1 + r_2 T)\frac{I}{A}, \\
    U_{act} & = s \log \left( \left(t_1 + \frac{t_2}{T} + \frac{t_3}{T^2} \right) \frac{I}{A} + 1 \right).
    \end{align}
\end{subequations}
$A$ is the electrode area of a single cell and $T$ is the temperature of the electrodes. The parameters $r_1$ and $r_2$ represent the ohmic resistance, while $s$ and $t_1,t_2,t_3$ are the activation overvoltage coefficients. Table \ref{tab:paramters_electrolyzer} reports the parameter estimates. The total power in a stack of $n_c$ electrolytic cells connected in series is 
\begin{equation} \label{eq:el3}
    P_{el} = I \, n_c \, U_{cell}.
\end{equation}
The total production flow rate of hydrogen of an electrolyzer with $n_c$ cells connected in series is a function of the current, according to
\begin{equation}
    f_{H_2} = n_c \left( \eta_F \frac{I}{z F} \right).
\end{equation}
$\eta_F$ is called Faraday efficiency (or current efficiency). This is experimentally defined as the ratio of the ideal hydrogen production rate to the actual hydrogen output from the electrolyzer. This is due to the fact that not all the electrons contribute to the water-splitting reaction. From the stoichiometry of the reaction, we notice that
\begin{equation}
    f_{H_2 O} = f_{H_2} = 2 f_{O_2}.
\end{equation}
The following model for the Faraday efficiency is proposed \cite{ULLEBERG200321}
\begin{equation} \label{eq:el4}
    \eta_F = \frac{(I/A)^2}{f_1 + (I/A)^2} f_2.
\end{equation}
$f_1$ and $f_2$ are parameters to be estimated empirically and they depend on the temperature. Table \ref{tab:paramters_electrolyzer} presents the estimated values for $T=80\, {}^\circ \text{C}$. 

\textbf{Remark.} The overvoltage $(U_{cell}-U_{tn})$ is responsible for generating heat in the electrolyzer (this is equal to $Q_{gen} = n_C (U_{cell}-U_{tn})I$). Therefore, a cooling system is generally deployed to keep the temperature constant. The exchanged heat could be stored or used, e.g., for district heating. We assume that the temperature in the electrolyzer is always kept constant.

\subsection{Hydrogen tank}
We model the hydrogen (and oxygen) tank considering the mass balance equation
\begin{equation} \label{eq:H2}
    \dot n_{H_2}^{stored} = f_{H_2}^{in} - f_{H_2}^{out}.
\end{equation}
The pressure in the tank can be obtained using a cubic equation of state (EOS), e.g. the Peng-Robinson EOS.

\begin{table}[tb]
    \centering
    \caption{Electrochemical model parameters values.}
    \begin{tabular}{|c|l l|}
        \hline
        \textbf{Parameter} & \textbf{Value} & \textbf{Unit}\\
        \hline
        \multicolumn{3}{|c|}{\textbf{U-I model \cite{SAKAS20224328}}}\\
        \hline
         $r_1$ & 0.8 & $\Omega\,\text{cm}^2$ \\
        $r_2$ & -0.00763 & $\Omega\,\text{cm}^2 \, {}^\circ\text{C}^{-1}$ \\
        $s$ & 0.1795 & $\text{V}$ \\
        $t_1$ & 20 &$ \text{cm}^2\, \text{A}^{-1}$ \\
        $t_2$ & 0.1 & $\text{cm}^2\, {}^\circ \text{C} \,\text{A}^{-1}$ \\
        $t_3$ & 3.5 $\times 10^5$ & $\text{cm}^2\,{}^\circ \text{C}^2 \,\text{A}^{-1}$ \\
        \hline
        \multicolumn{3}{|c|}{\textbf{Faraday efficiency model \cite{ULLEBERG200321}}}\\
        \hline
        $f_1 (T=80\, {}^\circ \text{C})$ & 250 & $\text{mA}^2 \text{cm}^{-4}$  \\
        $f_2 (T=80\, {}^\circ \text{C})$ & 0.980 & - \\
        \hline
    \end{tabular}
    \label{tab:paramters_electrolyzer}
\end{table}

\section{Optimal control} \label{sec:optimal_control}

We want to control the renewable energy system shown in Fig \ref{fig:diag}. The goal is to satisfy electricity demand from the grid and to maximize the profit generated from selling excess electricity, heat, and hydrogen. We also allow re-conversion of heat and hydrogen to electricity. This is achieved using a fuel cell for the hydrogen and a steam turbine for the heat coming from the thermal storage. Electricity can also be bought from the market and stored in the battery, when needed. The setup of the system is described in Fig. \ref{fig:diag}. We approximate the power generated by a generic fuel cell with $\eta_{fc}$ efficiency and hydrogen flow rate $d_{H_2}(t)$:
\begin{equation}\label{eq:fuel_cell}
        d_{fc}(t) = 141800 \cdot \eta_{fc} \cdot MM_{H_2} \cdot d_{H_2}(t).
\end{equation}
$MM_{H_2}$ is the molar mass of hydrogen.
We introduce auxiliary variables to simplify the formulation of the optimal control problem. 
Given the total produced power, $P_{tot}(t)$, by solar PVs and wind turbines and the electricity demand, $d_{el}(t)$, we define the excess power, $P_{sto}(t)=P_{tot} (t) - d_{el}(t)$ (when $P_{tot} (t) > d_{el}(t)$), which is to be stored. We define the variable $d_{sto}(t)$ as the power demand from the storage systems, when the produced power is insufficient ($d_{el}(t) > P_{tot}(t)$). Using these variables, the optimization problem only controls the energy surplus and deficit, while what is produced until surplus goes directly to the grid without passing through the storage. We define the following economic nonlinear optimal control problem
\begin{subequations} \label{eq:optimal_control}
\allowdisplaybreaks
\begin{align}
\max_{u(t),x(t)} \quad & \phi = \phi([x(t), u(t)]_{t_0}^{t_f})\\
\text{s.t.} \quad &  x(t_0) = \hat x_0, \\
& P_{sto}(t) = P_{b,in}(t)+P_{t,in}(t)+ P_{el}(t), \label{eq:constraint1}\\
& d_{sto}(t) = d_{b}(t) + d_{t}(t) +  d_{fc}(t), \label{eq:constraint2}\\
&\dot E_{b}(t) = -\alpha_b E_b(t) + \eta_{b,in} \left( P_{b,in}(t) + P_{b,bo}(t) \right) \nonumber \\ & \quad \quad \quad \quad - \frac{1}{\eta_{b,out}} \left( P_{b,out}(t) + d_{b}(t) \right), \\
&	m \, c_{P}(T) \dot T(t) = \eta_{t,in} P_{t,in}(t) - (UA)_{t} (T(t)-T_a) \nonumber \\ & \quad \quad \quad \quad - Q_{t,out}(t) - \frac{1}{\eta_{steam}} d_{t}(t), \\
& \text{Electrolyzer model \eqref{eq:el1}, \eqref{eq:el2}, \eqref{eq:el3}, \eqref{eq:el4}} \nonumber, \\
& \text{Fuel cell model \eqref{eq:fuel_cell}}, \nonumber \\
& \dot n_{H_2}^{stored}(t) = f_{H_2}^{in}(t) - f_{H_2}^{out}(t) - d_{H_2}(t), \\
& x_{min} \leq x(t) \leq x_{max}, \\
& u_{min} \leq u(t) \leq u_{max}. 
\end{align}
\end{subequations}
The objective function is
\begin{equation}
\begin{split}
     \phi & =  \int_{t_0}^{t_f} c^e(t) P_{b,out}(t) \, dt - \int_{t_0}^{t_f} c^e(t) P_{b,bo}(t) \, dt +\\
     & \int_{t_0}^{t_f} c^h(t) Q_{t,out}(t) \, dt + MM_{H_2} \int_{t_0}^{t_f} c^{H_2}(t) f_{H_2}^{out}(t) \, dt \\
     & + \text{Profit to-go}.
\end{split}
\end{equation}
$c^e(t), c^h(t)$ and $c^{H_2}(t)$ are the cost functions of electricity, heat and hydrogen over time, respectively. The subscript $out$ indicates outgoing flows to the markets (energy sold), while $in$ indicates incoming flows to the storage devices from wind and solar. Constraint \eqref{eq:constraint1} dispatches the energy surplus to the storage devices. Constraint \eqref{eq:constraint2} ensures that power is withdrawn from the storage devices when the production is not sufficient. This constraint can be formulated as a soft constraint, to ensure feasibility of the problem even when there is not enough energy in the system. Notice that the nonlinearity is due to the dynamics in the electrolyzer.

The states consist of the dynamics in the ESSs:
\begin{equation}
    x(t) = \Big[ E_b(t), T(t), n_{H_2}^{stored}(t) \Big].
\end{equation}
The vector of manipulated variables is
\begin{equation}
\begin{split}
    u(t) = & \Big[ P_{b,in}(t), P_{t,in}(t),P_{el}(t), P_{b,out}(t), P_{b,bo}(t), \\ &  Q_{t,out}(t), f_{H_2}^{out}(t),I(t) \Big].
\end{split}
\end{equation}
The disturbances of the system are
\begin{equation}
    d(t) = \Big[ c^e(t), c^h(t), c^{H_2}(t), P_{sto}(t), d_{sto}(t) \Big].
\end{equation}

\section{Simulation results} \label{sec:simulation}
\begin{figure}[tb]
    \centering
    \includegraphics[width=0.48\textwidth]{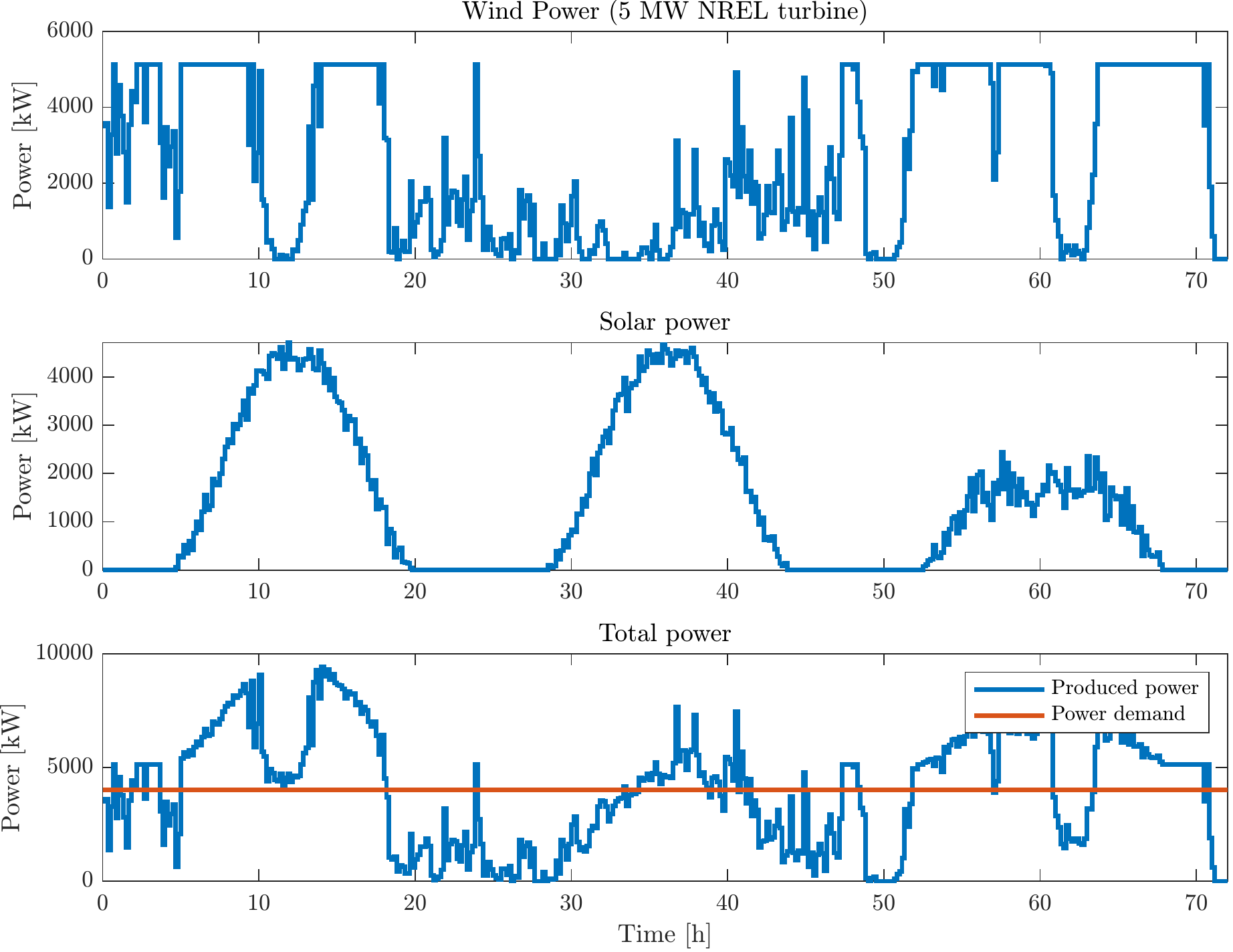}
    \caption{Simulation of the amount of power produced from wind and solar over three days. Only the mean wind speed is considered.}
    \label{fig:3days}
\end{figure}

\begin{figure}[tb]
    \centering
    \includegraphics[width=0.49\textwidth]{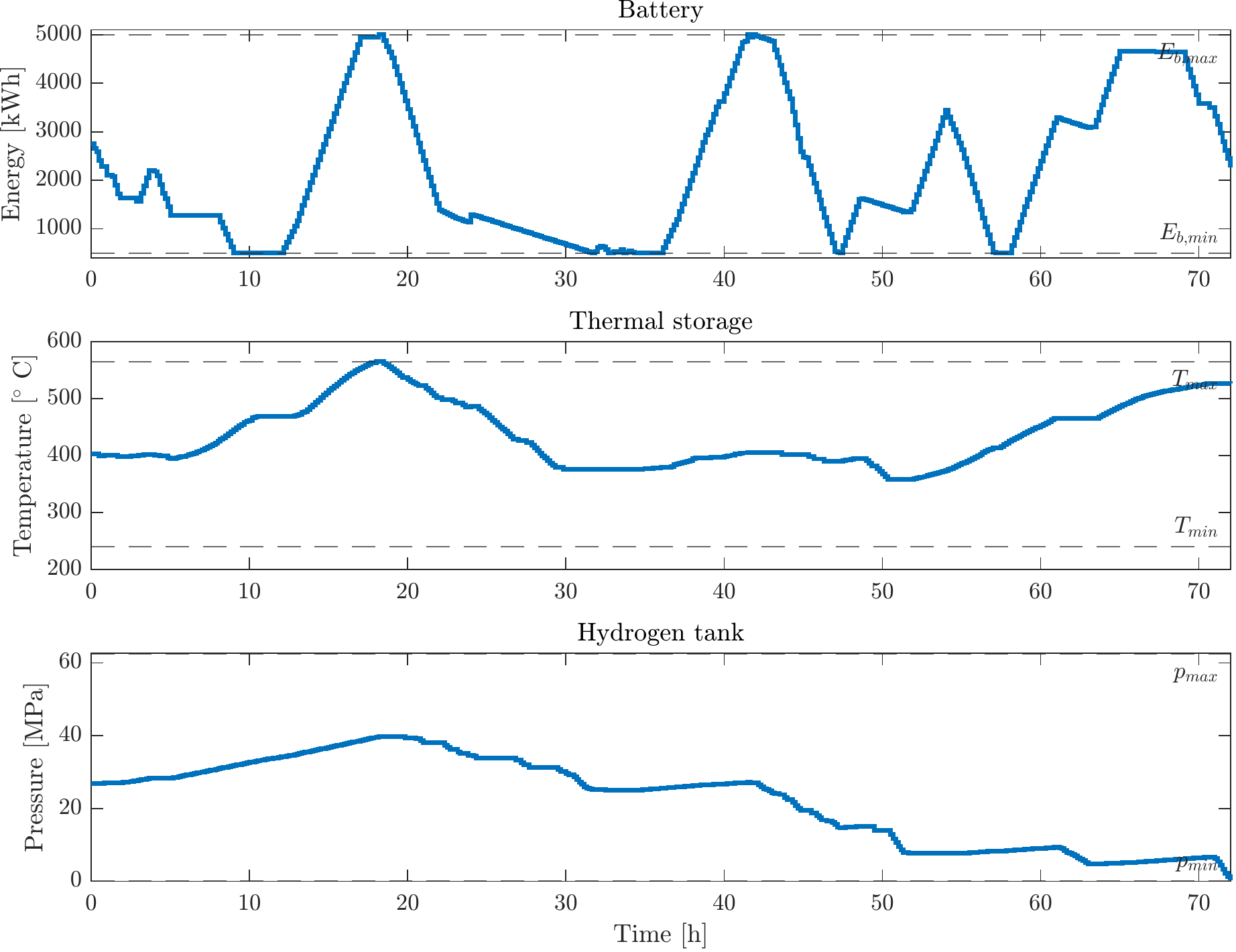}
    \caption{Solution of the optimal control problem: states $x(t)$.}
    \label{fig:states}
\end{figure}

\begin{figure}[tb]
    \centering
    \includegraphics[width=0.49\textwidth]{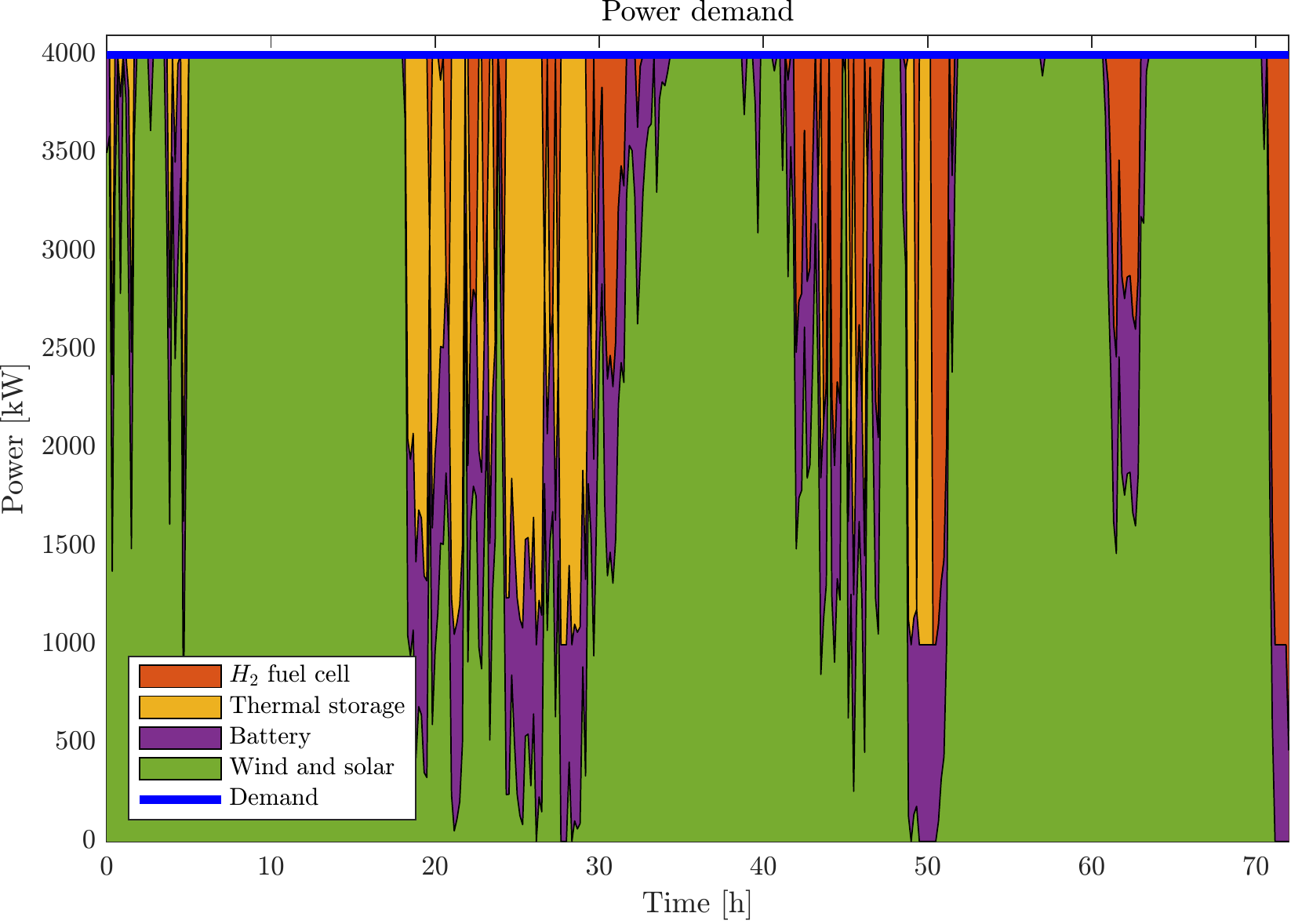}
    \caption{Solution of the optimal control problem: power balance of the energy system.}
    \label{fig:power_demand}
\end{figure}

We present a 3 days simulation. A constant power demand $d_{el}(t)$ of 4 MW is assumed. The configuration of the system is the following: one NREL 5 MW turbine, a 6 MW solar PV park, 50 MWh thermal storage, 5 MWh battery, 2.4 MW electrolyzer, and a 1000 kg (30 m$^3$) $\simeq$ 39 MWh hydrogen tank. The time series for the prices are randomly generated from normal distributions. We remark that modeling energy prices is outside of the scope of this paper. Our main purpose is to demonstrate optimal control of renewable energy systems with storage.
Fig. \ref{fig:3days} shows the simulated power production (solar, wind and total).
The optimal control problem \eqref{eq:optimal_control} is discretized, using 10 min as sampling time and 1 h as control interval. The problem is set up using CasADi \cite{Andersson2019} and solved using IPOPT. Fig. \ref{fig:states} presents the states $x(t)$ of the system. 
Fig. \ref{fig:power_demand} displays how the electricity demand is satisfied by the control program, i.e. the power distribution. The program controls the power flows such that the demand is always met, optimally withdrawing energy from the storage systems.
\section{Conclusions}
\label{sec:Conclusions}
In this paper, we present models of energy production from wind and solar and energy storage via a battery, thermal storage, and electrolysis. Such an energy system can be operated optimally using the optimal control program that we formulate in section \ref{sec:optimal_control}. This program can be used as the regulator in an MPC algorithm.
The program smartly charges the storage systems before a low power production is expected, and discharges them to satisfy the demand. The excess power is sold when possible, and when prices are high.
We demonstrate that the power demand can be successfully satisfied, while maximizing the economic revenue. Incorporating multiple and different ESSs is beneficial to renewable energy systems, providing flexibility and reliability. Moreover, they increase the economic value of the system, by eliminating or minimizing energy curtailment, storing the surplus of energy. 

\bibliographystyle{ieeetr}
\bibliography{Bibliography}

\end{document}